\documentstyle{article}
\title{Projective normality of abelian surfaces of type $(1,2d)$.}
\author{Luis Fuentes Garc{\'{\i}}a \thanks{Supported by EAGER.}}
\date{}

\newtheorem{teo}{Theorem}[section]
\newtheorem{defin}[teo]{Definition}
\newtheorem{prop}[teo]{Proposition}

\newtheorem{lemma}[teo]{Lemma}

\newtheorem{rem}[teo]{Remark}

\def\E{{\cal E}}
\def\Te{{\cal O}}

\def\L{{\cal L}}
\def\M{{\cal M}}
\def\LY{{\cal L}_1}
\def\LYY{{\cal L}_2}

\def\P{{\bf P}}

\def\H{{\bf H}}
\def\Z{{\bf Z}}

\def\qed{\hspace{\fill}$\rule{2mm}{2mm}$}

\def\lrw{{\longrightarrow}}

\def\dual{{\star}}

\begin{document}

\maketitle

{\bf Abstract:} We show that an abelian surface embedded in $\P^N$
by a very ample line bundle $\L$ of type $(1,2d)$ is projectively
normal if and only if $d\geq 4$. This completes the study of the
projective normality of abelian surfaces embedded by complete
linear systems. 

\vspace{0.1cm}

\section{Introduction}\label{intro}

Let $A$ be an abelian surface. Let $\L$ be an ample line bundle on
$A$ of type $(n_1,n_2)$. It induces a rational map
$\phi_{\L}:A\lrw \P^{n_1n_2-1}$. We want to study when this map is
a projectively normal embedding. The known results are the
following:

\begin{enumerate}

\item If $n_1\geq 3$ then $\phi_{\L}$ is a projectively normal
embedding (see \cite{ko}, \cite{bila}).

\item If $n_1=2$ then $\phi_{\L}$ is a projectively normal
embedding if and only if no point of $K(\L)$ is a base point of
$\L'$, where $\L=\L'^2$ (see \cite{oh}, \cite{bila}).

\item If $n_1=1$ then $\L$ is a primitive bundle of type
$(1,n_2)$. In this case:

\begin {enumerate}

\item If $n_2=7,9,11$ or $n_2\geq 13$ then $\phi_{\L}$ is a
projectively normal embedding if and only if $\L$ is very ample
(see \cite{la}).

\item If $n_2\geq 7$ and $A$ is generic (in particular
$NS(A)\simeq \Z$) then $\phi_{\L}$ is a projectively normal
embedding (see \cite{iy1}).

\item If $n_2> 8$ and $A$ is not isogenous to product of elliptic curves then $\phi_{\L}$ is a projectively normal
embedding (see \cite{iy2}).

\end{enumerate}

\end{enumerate}
Note that if $n_1n_2<7$ the embedding can never be projectively
normal, because the dimension of $Sym^2H^0(\L)$ is less than that
of $H^0(\L^2)$. Thus the open cases are when the line bundle $\L$
is of type $(1,8)$, $(1,10)$ or $(1,12)$. In this paper we study
these cases. The main result is the following:

\begin{teo}
Let $A$ be an abelian surface and let $\L$ be a line bundle on $A$
of type $(1,n)$. The induced map $\phi_{\L}:A\lrw \P^{n-1}$ is a
projectively normal embedding if and only if $\L$ is very ample
and $n\geq 7$.
\end{teo}

It is well known that a necessary condition for $\phi_{\L}$ to be
a projectively normal embedding is the very ampleness of $\L$.
Furthermore, by Proposition $2.3$, \cite{iy2} it is sufficient to
analyze the $2$-normality of the map $\phi_{\L}$, that is, study
the quadrics containing the image of the abelian surface.

We will work with a very ample line bundle $\L$ of type $(1,2d)$.
Taking a nontrivial $2$-torsion point $x\in K(\L)$, we construct
an involution on $A$ that extends to $\P^{2d-1}$. From this, there
are two disjoint spaces of fixed points in $\P^{2d-1}$ of
complementary dimension. Moreover there is an induced involution
in the space of quadrics of $\P^{2d-1}$. This space decomposes
into two subspaces of invariant quadrics. We call one of them the
space of base quadrics: these are quadrics containing the two
spaces of fixed points. The other one is the space of harmonic
quadrics: they have the property that the spaces of fixed points are polar
spaces with respect to them. We will see how this decomposition is
related to the decomposition into Heisenberg modules.

On the other hand, we will consider the scroll $R$ obtained by
joining with lines the points of $\phi_{\L}(A)$ related by the
involution. We will see that $R$ and $\phi_{\L}(A)$ are contained
in the same number of independent base quadrics. We will modify
the arguments of R. Lazarsfeld in \cite{la} to bound this number.
This will allow us to solve the case $(1,10)$. Moreover, we will obtain that the abelian surface embedded by a
very ample line bundle of type $(1,6)$ is never contained
in quadrics.

Finally, we will work with the space of harmonic quadrics
containing the abelian surface to complete the cases of the
polarizations of type $(1,8)$ and $(1,12)$.

{\bf Acknowledgement} I thank K. Hulek and H.-Ch. von Bothmer for
their useful comments and suggestions. I am also grateful to the
Institut f\"{u}r Mathematik of Hannover for its hospitality and to the project EAGER for financial support.

\section{Preliminaries.}\label{quadricsinvolution}

First, let us recall some general facts about involutions and
double covers of smooth varieties:

\begin{defin}
Let $X,Y$ be two smooth varieties. A map $\pi:X\lrw Y$ is called
a double cover of $Y$ if it is a finite map of degree $2$.
\end{defin}

\begin{teo}\label{doublecover}
Let $\pi:X\lrw Y$ be a double cover. Then:
\begin{enumerate}

\item $\pi_*\Te_X\simeq \Te_Y\oplus \M$ where $\M$ is a line
bundle on $Y$.

\item $\M^2\simeq \Te_Y(-B)$ where $B$ is the branch divisor of
$\pi$.
\end{enumerate}
\end{teo}
{\bf Proof:} See \cite{pe}. \qed

\begin{defin}
We say that a line bundle $\L$ on $X$ is invariant by the double
cover $\pi:X\lrw Y$ when $\L\simeq \pi^*\LY$, where $\LY$ is a
line bundle on $Y$. Equivalently, $\L$ is invariant by $\pi$ if
and only if $\L\simeq \mu^*(\L)$, where $\mu$ is the involution of
$X$ induced by $\pi$.
\end{defin}

From now on we will work with a double cover $\pi:X\lrw Y$ and an
invariant line bundle $\L$ on X. The map $\mu:X\lrw X$ will be the
involution induced by $\pi$.

Let us consider the rank $2$ vector bundle $\E=\pi_* \L$ and let
$S=\P(\E)$ be the corresponding $\P^1$ bundle on $Y$. Moreover,
take $\E'=\pi_*(\L^2)$. From the definition of invariant line
bundle, the projection formula and Theorem \ref{doublecover} it follows that:
$$
\begin{array}{l}
{\E=\pi_*\L=\pi_*\pi^*\LY=\LY\oplus \LYY}\\
{\E'=\pi_*(\L^2)=\pi_*\pi^*(\LY^2)=\LY^2\oplus (\LY\otimes \LYY)}\\
\end{array}
$$
where $\LYY= \LY\otimes \M$.

From this, we have canonical isomorphisms of vector spaces:
$$
\begin{array}{l}
{H^0(X,\L)\simeq H^0(S,\Te_S(1))\simeq H^0(Y,\E)\simeq
H^0(Y,\LY)\oplus H^0(Y,\LYY)}\\
{H^0(X,\L^2)\simeq H^0(Y,\E')\simeq H^0(Y,\LY^2)\oplus
H^0(Y,\LY\otimes \LYY)}\\
{H^0(S,\Te_S(2))\simeq H^0(Y,S^2\E)\simeq H^0(Y,\LY^2)\oplus H^0(Y,\LY\otimes
\LYY)\oplus H^0(Y,\LYY^2).}\\
\end{array}
$$

Note that $H^0(Y,\LY)$ and $H^0(Y,\LYY)$ correspond to the $\pm 1$-eigenspaces of $H^0(X,\L)$ with respect to the involution $\mu$.

In order two study the quadrics containing the projective images
of $X$ and $S$, we consider the maps:
$$
\begin{array}{l}
{\alpha:Sym^2(H^0(X,\L))\lrw H^0(X,\L^2)}\\
{\beta:Sym^2(H^0(S,\Te_S(1)))\lrw H^0(S,\Te_S(2)).}\\
\end{array}
$$
Using the previous isomorphisms, we see that they decompose into the
following way:
$$
\begin{array}{l}
{\alpha\simeq \alpha_h\oplus \alpha_b }\\
{\beta\simeq \beta_h\oplus \beta_b}\\
\end{array}
$$
with,
$$
\begin{array}{l}
{\alpha_h:Sym^2(H^0(Y,\LY))\oplus Sym^2(H^0(Y,\LYY))\lrw
H^0(Y,\LY^2)}\\
{\alpha_b:H^0(Y,\LY)\otimes H^0(Y,\LYY)\lrw H^0(Y,\LY\otimes \LYY)}\\
\end{array}
$$
and
$$
\begin{array}{l}
{\beta_h:Sym^2(H^0(Y,\LY))\oplus Sym^2(H^0(Y,\LYY))\lrw
H^0(Y,\LY^2)\oplus H^0(Y,\LYY^2)}\\
{\beta_b:H^0(Y,\LY)\otimes H^0(Y,\LYY)\lrw H^0(Y,\LY\otimes \LYY).}\\
\end{array}
$$
Let us remark that $\alpha_h$ and $\beta_h$ correspond to the restriction of $\alpha$ and $\beta$ to the $1$-eigenspaces with respect to the involutions induced by $\mu$ on $Sym^2(H^0(Y,\LY^2))$ and $Sym^2(H^0(S,\Te_S(1)))$. On the other hand, $\alpha_b$ and $\beta_b$ correspond to the restriction to the $-1$-eigenspaces. To understand the map $\alpha_h$, note that, by Theorem
\ref{doublecover}, there is a natural inclusion of linear systems
$$
|\LYY^2|\subset |\LY^2|
$$
given by the addition of the branch divisor $B$.

There is a nice geometrical interpretation of these maps. Suppose
that the line bundle $\L$ is very ample. Then the linear system $|\Te_S(1)|$
separates the points related by the involution on each ruling of $S$. Therefore, $|\Te_S(1)|$ is base point free linear system. Furthermore, the $\P^1$-bundle $S$ has two natural sections given by the epimorphisms $\LY\oplus \LYY\lrw \LY$ and $\LY\oplus \LYY\lrw \LYY$. The linear system $|\Te_S(1)|$ restricted to them correspond to the linear systems $|\LY|$ and $|\LYY|$ on $Y$, so they are base point free too.

We have an embedding $X\lrw \P^N$ where
$\P^N=\P(H^0(X,\L)^\dual)$. We will identify $X$ with its image.
On the other hand, we have a map $S\lrw \P(H^0(X,\L)^\dual)$.
This is not necessarily birational, but anyway, the image is a
scroll. We will denote it by $R$. The involution $\mu$ on $X$
extends to $\P^N$. There are two subspaces of fixed points:
$\P(H^0(Y,\LY)^\dual)$ and $\P(H^0(Y,\LYY)^\dual)$.

Moreover, we have an induced involution in the space of quadrics
of $\P^N$. This space decomposes into two subspaces of invariant
quadrics corresponding to the $\pm 1$-eigenspaces:
$$
H^0(\Te_{P^N}(2))\simeq H^0(\Te_{P^N}(2))_h\oplus
H^0(\Te_{P^N}(2))_b
$$
where
$$
\begin{array}{l}
{H^0(\Te_{P^N}(2))_h\simeq Sym^2(H^0(Y,\LY))\oplus
Sym^2(H^0(Y,\LYY))}\\
{H^0(\Te_{P^N}(2))_b\simeq H^0(Y,\LY)\otimes H^0(Y,\LYY).}\\
\end{array}
$$

\begin{defin}
With the previous notation, we call $H^0(\Te_{P^N}(2))_h$ the
space of harmonic quadrics of $\P^N$ and we call $H^0(\Te_{P^N}(2))_h$
the space of base quadrics of $\P^N$.
\end{defin}

\begin{rem}
The base quadrics correspond to quadrics containing the spaces of
fixed points of $\P^N$; on the other hand, these spaces are polar
spaces respect to the harmonic quadrics. \qed
\end{rem}

Now, the spaces $ker(\alpha_b)$ and $ker(\beta_b)$ corresponds to
the spaces of base quadrics containing $X$ and $R$ respectively
and the spaces $ker(\alpha_h)$ and $ker(\beta_h)$ corresponds to
the spaces of harmonic quadrics containing $X$ and $R$
respectively.

Note, that $ker(\alpha_b)\simeq \ker(\beta_b)$. This corresponds
to the following geometrical fact: a quadric containing $X$ and
a space of fixed points meets each generator of $R$ in three
points, so it contains the scroll. Thus, the spaces of base
quadrics containing $X$ and $R$ are the same.

We finish this section by proving an useful result to bound the
number of independent quadrics containing a projective variety.

\begin{lemma}\label{quadrics}
Let $X\subset \P^N$ be an irreducible non-degenerated
$d$-dimensional variety of degree $m$. Let $k$ be the number of
independent quadrics containing $X$. Then:
\begin{enumerate}
\item If $m>2(N-d)$ then $k\leq \left(^{N-d+2}_{\quad\;
2}\right)-2(N-d)-1$ \item If $m\leq 2(N-d)$ then $k\leq
\left(^{N-d+2}_{\quad\; 2}\right)-m.$
\end{enumerate}
\end{lemma}
{\bf Proof:} Let $\P^{N-d}$ be a generic $(N-d)$-dimensional space
in $\P^N$. It cuts $X$ in a set $M$ of $m$ points in general
position. They are contained in at least the same number of
independent quadrics than $X$. Moreover, it is well known that
$2r+1$ points in $\P^r$ in general position impose independent
conditions to the quadrics. From this if $m> 2(N-d)$ we have
$k\leq h^0(\Te_{P^{N-d}}(2))-2(N-d)-1$ and if $m\leq 2(N-d)$ then
$k\leq h^0(\Te_{P^{N-d}}(2))-m$. \qed

\section{The Heisenberg decomposition of the space of
quadrics.}\label{heisenberg}

Let $V$ be a $n$-dimensional vector space. Let us take coordinates
$\{x_0,\ldots,x_{n-1}\}$ in $V$. Let us consider the
Schr\"{o}dinger representation of the Heisenberg group $\H_{n}$.
We take generators $\sigma,\tau\in GL(V)$, defined by:
$$
\begin{array}{l}
{\sigma(x_j)=x_{j-1}}\\
{\tau(x_j)=\varepsilon^j x_j \mbox{, with }\varepsilon=e^{2\pi \imath /n}}.\\
\end{array}
$$

The Heisenberg group acts on $Sym^2V$. This corresponds to the
space of quadrics of $\P(V)$. We will denote it by $W$. The
dimension of $W$ is $n(n+1)/2$. We will describe explicitly the
Heisenberg decomposition of $W$. We distinguish two cases
depending of the parity of $n$:

\begin{enumerate}

\item Suppose that $n$ is odd.

In this case $W$ decomposes into $\frac{n+1}{2}$ $n$-dimensional
irreducible $\H_{n}$-modules:
$$
W\simeq \bigoplus W_l
$$
where
$$
W_l=\langle x_ix_j/i=j+l; i,j\in \Z_{n} \rangle
$$
and $l=0,\ldots,\frac{n+1}{2}-1$. Because all of them are
isomorphic as $\H_{n}$-modules, we will write $W_l\simeq V'$.

\item Suppose that $n$ is even. Let us take $n=2d$.

By using the Schr\"{o}dinger representation, we find the following
decomposition of $W$ in four $\H_{n}$-modules, $W\simeq
W_0^+\oplus W_0^-\oplus W_1^+\oplus W_1^-$, with,
$$
\begin{array}{l}
{W_0^+=\langle x_ix_j+x_{i+d}x_{j+d}/i+j\equiv_2 0; i,j\in \Z_{2d} \rangle}\\
{W_0^-=\langle x_ix_j-x_{i+d}x_{j+d}/i+j\equiv_2 0; i,j\in \Z_{2d} \rangle}\\
{W_1^+=\langle x_ix_j+x_{i+d}x_{j+d}/i+j\equiv_2 1; i,j\in \Z_{2d} \rangle}\\
{W_1^-=\langle x_ix_j-x_{i+d}x_{j+d}/i+j\equiv_2 1; i,j\in \Z_{2d} \rangle}.\\
\end{array}
$$
They are not irreducible, but each one of them decomposes into
isomorphic $d$-dimensional irreducible $\H_{n}$-modules:
$$
{W_l^\pm\simeq \bigoplus W_{l,m}^\pm}
$$
where
$$
{W_{l,m}^\pm=\langle x_ix_j\pm x_{i+d}x_{j+d}/i+j\equiv_2 l; i=j+m; i,j\in
\Z_{2d} \rangle}\\
$$
and $l=0,1$; $m=0,\ldots,d$; $l\equiv_2 m$. If we fix $l$ and the
sign $\pm$, the spaces $W_{l,m}^\pm$ are isomorphic as
$\H_{n}$-modules. We will write $W_{l,m}^\pm\simeq V_l^\pm$.

\end{enumerate}

We summarize these results in the following theorem.

\begin{teo}\label{decomposition}
Let $V$ be a vector space of dimension $n$. Let $\H_{n}$ be the
$n$-dimensional Heisenberg group acting on $V$. The space
$W=Sym^2V$ has the following decomposition into $\H_n$-modules:
\begin{enumerate}

\item If $n$ is odd then
$$W\simeq \frac{n+1}{2} V'$$
where $V'$ is a $n$-dimensional irreducible $\H_n$-module.

\item If $n=2d$ is even then:

\begin{enumerate}

\item If $d$ is odd,
$$W\simeq \frac{d+1}{2} V_0^+\oplus
\frac{d+1}{2} V_0^-\oplus \frac{d+1}{2} V_1^+\oplus \frac{d-1}{2}
V_1^-.$$

\item If $d$ is even,
$$
W\simeq \frac{d+2}{2} V_0^+\oplus \frac{d}{2} V_0^-\oplus
\frac{d}{2} V_1^+\oplus \frac{d}{2} V_1^-.
$$

\end{enumerate}

The spaces $V_k^\pm$ are $d$-dimensional $\H_n$-modules and no two
of them are isomorphic. \qed

\end{enumerate}
\end{teo}

\begin{rem}\label{four}
Let us suppose that $n=2d$. In this case $\sigma^d$ and $\tau^d$
are elements of order two in $\H_{n}$ and they define involutions
$\mu_{\sigma}$, $\mu_{\tau}$ in $\P(V)$. Moreover,
$\sigma^d\tau^d$ is not necessarily of order two, but
$(\sigma^d\tau^d)^2=(-1)^d Id_V$ so it defines an involution
$\mu_{\sigma\tau}$ in $\P(V)$.

The space of quadrics $W$ decomposes with respect to these
involutions into two subspaces of invariant quadrics. With the
previous notation we have:

\begin{enumerate}

\item With respect to the involution $\mu_{\sigma}$, $W=\langle
W_0^+,W_1^+\rangle\oplus \langle W_0^-,W_1^- \rangle$ of dimensions $d(d+1)$ and $d^2$ respectively.

\item With respect to the involution $\mu_{\tau}$, $W=\langle
W_0^+,W_0^-\rangle\oplus \langle W_1^+,W_1^- \rangle$ of dimensions $d(d+1)$ and $d^2$ respectively.

\item With respect to the involution $\mu_{\sigma\tau}$,
$W=\langle W_0^+,W_1^-\rangle\oplus \langle W_0^-,W_1^+ \rangle$ of dimensions $d^2$ and $d(d+1)$ respectively if $d$ is odd and $d(d+1)$ and $d^2$ respectively if $d$ is even.

\end{enumerate}

Moreover, we know that the dimension of the spaces of harmonics quadrics and base quadrics is $d(d+1)$ and $d^2$ respeectively. If we compare these dimensions with the
dimensions of the invariant spaces in each case we see the
following: when $d$ is odd the space of base quadrics must contain
$W_1^-$; when $d$ is even the space of harmonic quadrics must
contain $W_0^+$. \qed
\end{rem}

\section{Quadrics containing abelian varieties of type
$(1,\ldots,1,2d)$.}\label{abelianvarieties}

Let $A$ be an abelian variety of dimension $g$ and let $\L$ be an
ample line bundle of type $(1,\ldots,1,2d)$. Let $K(\L)$ be the
kernel of the isogeny $\lambda_{\L}:A\lrw \hat{A}\simeq Pic^0(A)$
determined by the line bundle $\L$. It is well known that
$K(\L)\simeq \Z_{2d}\times \Z_{2d}$. The natural extension
of $K(\L)$ by $\langle \varepsilon \rangle\subset {\bf C}^*$, with $\varepsilon=e^{2\pi \imath /(2d)}$, is the Heisenberg group $\H_{2d}$ of dimension $2d$.

Suppose that $|\L|$ is very ample. It induces an embedding
$\phi_{\L}:A\lrw \P^{2d-1}=\P(V)$, with $V=H^0(A,\L)^\dual$. We
will identify $A$ with its image. Let $I$ be the space of quadrics
containing $A$. This is an $\H_{2d}$-submodule of the space of
quadrics of $\P^{2d-1}$. By Theorem \ref{decomposition}, the space
$I$ decomposes into the following $\H_{2d}$-modules:
$$
I=I_0^+\oplus  I_0^-\oplus I_1^+\oplus I_1^-
$$
where $I_j^\pm\subset W_j^\pm$.

On the other hand, let $\{x_1,x_2,x_3\}$ be the three nontrivial $2$-torsion points
in $K(\L)$. Let $x=x_i$ be one of them. We can consider the double
cover $\pi:A\lrw \bar{A}\simeq A/x$. $\L$ is invariant by $\pi$.
We are in the situation described in section
\ref{quadricsinvolution}. Now, $\bar{A}$ is an abelian variety;
$\L=\pi^*\LY$, where $\LY$ is a line bundle in $\bar{A}$ of type
$(1,\ldots,1,d)$; $\LYY=\LY\otimes \M$, where $\M\in
Pic^0(\bar{A})$ verifies $M^2\simeq \Te_{\bar{A}}$. Joining by
lines the points of $A$ related by the involution we obtain a
scroll $R\subset \P^{2d-1}$. $R$ is the image of the $\P^1$ bundle
$S=\P(\LY\oplus \LYY)$ by the tautological linear system
$|\Te_S(1)|$ (see \cite{cihu} for a detailed construction).
Moreover, $I$ decomposes into a space $I_b$ of base quadrics and a
space $I_h$ of harmonic quadrics:
$$
\begin{array}{l}
{I_b=ker(\alpha_b:H^0(\LY)\oplus H^0(\LYY)\lrw H^0(\LY^2))}\\
{I_h=ker(\alpha_h:Sym^2(H^0(\LY))\oplus Sym^2(H^0(\LYY))\lrw H^0(\LY\otimes
\LYY)).}\\
\end{array}
$$

Note that this decomposition depends on the choice of the $2$-torsion point $x\in K(\L)$. In fact, the three involutions defined on $\P(V)$ by $x_1,x_2,x_3$ correspond to the three involutions described in Remark \ref{four}. In this way, we know how this decomposition is related to the Heisenberg decomposition:

\begin{prop}\label{bh}
The spaces $I_h$ and $I_b$ decompose into the sum of two
$H_{2d}$-modules. In particular:

\begin{enumerate}
\item If $d$ is odd, $I_b=I_1^-\oplus I_0^-$, $I_b=I_1^-\oplus I_1^+$ or $I_b=I_1^-\oplus I_0^+$ where each possibility occurs for exactly one of the $2$-torsion points $x_1$, $x_2$ and $x_3$.

\item If $d$ is even, $I_h=I_0^+\oplus I_1^+$, $I_h=I_0^+\oplus I_0^-$ or $I_h=I_0^+\oplus I_1^-$ where each possibility occurs for exactly one of the $2$-torsion points $x_1$, $x_2$ and $x_3$.
\qed
\end{enumerate}
\end{prop}

\section{Quadrics containing abelian
surfaces of type $(1,2d)$.}\label{abeliansurfaces}

Let $A$ be an abelian surface and let $\L$ be a very ample line
bundle of type $(1,n)$. We want to study the projective normality
of the embedding given by the line bundle $\L$. By Proposition
$2.3$, \cite{iy2} it is sufficient to analyze the normality with
respect to quadrics.

Note that, if $n<5$, $\L$ is never very ample. When $n=5,6$, $\L$
cannot be $2$-normal because the dimension of $Sym^2(H^0(\L))$ is
less than the dimension of $H^0(\L^2)$. However, one can ask if
the map
$$Sym^2(H^0(\L))\lrw H^0(\L^2)$$ has the expected behavior, that
is, if it is injective. This is equivalent to study whether there are
quadrics containing the abelian surface.

The case $n=5$ is immediate. $A$ is embedded as a surface of
degree $10$ in $\P^4$. By the Heisenberg decomposition of the
space of quadrics, if the ideal of quadrics containing $A$ is not
empty, it has at least $5$ independent quadrics, and this is not
possible.

When $n>5$, we use the idea of R. Lazarsfeld in \cite{la}. He
applies Lemma \ref{quadrics} to the abelian surface and the fact
that the space of quadrics decomposes into $d$-dimensional $H_{n}$
irreducible modules, with $d=n/m.c.d.(n,2)$. This provides
directly the projective normality when $n=7,9,11$ and $n\geq 13$.

However, when $n=6,8,10,12$ this argument does not imply
projective normality, but it reduces the possibilities for the
number of independent quadrics containing $A$. If we call this number $k$, we obtain:
$$
\begin{tabular}{|c|c|} \hline
  {Polarization} & {$k$} \\ \hline
  {$(1,6)$} & {$0,3$} \\
  {$(1,8)$} & {$4,8$} \\
  {$(1,10)$} & {$15,20$} \\
  {$(1,12)$} & {$30,36$} \\ \hline
\end{tabular}
$$

We will analyze these cases in detail. We work with very ample
line bundles $\L$ of type $(1,n)$ with $n=2d$. We will use the
construction developed in the previous sections.

 Let us consider the involution given by a nontrivial
$2$-torsion point in $K(\L)$. The line bundles $\LY$ and $\LYY$ on
the abelian surface $\bar{A}$ are of type $(1,d)$. We have seen that they are base-point-free so they define regular
maps $\phi_1: \bar{A}\lrw \P^{d-1}_1$ and $\phi_2: \bar{A}\lrw
\P^{d-1}_2$. In particular, when $d\geq 4$ they are birational embeddings or
$2:1$ maps onto elliptic scrolls (see \cite{bilast}; \cite{hula};
\cite{ra}). When $d=3$, they are $6:1$ coverings of $\P^2$.
Moreover, $\P^{d-1}_1$ and $\P^{d-1}_2$ are disjoint subspaces of
$\P^{n-1}$. The scroll $R$ meets each one of these subspaces
exactly in the image of the maps $\phi_1$ and $\phi_2$.

We will study the number of independent base quadrics and the
number of independent harmonic quadrics containing the abelian
surface. We will denote these numbers by $k_b$ and $k_h$
respectively. We will need the following result:

\begin{lemma}\label{tecnico}
Let $C_2$ (resp $C_1$) be a generic curve in $|\LYY|$ (resp.
$|\LY|$). If $d\geq 4$, then the restriction of $\phi_1$ (resp.
$\phi_2$) to $C_2$ (resp. $C_1$) is a birational map.
\end{lemma}
{\bf Proof:} If the map $\phi_1$ is birational the result is
immediate, because $\LY, \LYY$ differ by translation.

Suppose that $\phi_1$ (and then $\phi_2$) is not birational.
Because $\LY$ and $\LYY$ are base point free and $d\geq 4$, it is
known that $\phi_1, \phi_2$ are $2:1$ maps onto elliptic scrolls (see \cite{hula}).
In particular there are involutions $\sigma_1,\sigma_2$ of
$\bar{A}$ such that $\sigma_i^*\L_i= \L_i$ and $\phi_i=\phi_i\circ \sigma_i$.

Suppose that the restriction of $\phi_1$ to any divisor $C_2\in
|\LYY|$ is not a birational map. This means that $\sigma_1=\sigma_2$ and in particular,
$\sigma_1^*\LYY= \LYY$. But $\LYY= \LY\otimes \M$, so it
holds that $\sigma_1^* \M= \M$. The abelian surface $A$ is
given locally by an equation $z_j^2-1=0$ in the line bundle $\M=
\bar{A}$ (see \cite{pe}). Since $\M$ is invariant by $\sigma_1$,
this involution can be lifted to an involution $\sigma$ in $A$. From this:
$$
\sigma^*\L= \sigma^*\pi^*\LY = \pi^*\sigma_1^*\LY=
\pi^*\LY= \L.
$$
We see that $\L$ is invariant by $\sigma$. Moreover, because $\phi_i=\phi_i\circ \sigma_i$ for $i=1,2$, the map $\phi$ defined by $\L$ verifies $\phi=\phi\circ \sigma$. But this contradicts its very ampleness. \qed

\begin{rem}\label{tecnico2}
When $d=3$, the maps $\phi_1,\phi_2$ are $6:1$ coverings over
$\P^2$. Since $\LYY\not\simeq \LY$, the generic divisor of
$|\LYY|$ (resp. $|\LY|$) is mapped by $\phi_1$ (resp. $\phi_2$)
to a nondegenerate curve of $\P^2$. \qed
\end{rem}

\subsection{Abelian surfaces of type $(1,6)$ and $(1,10)$.}

First, we will use the fact that there are the same number of
independent base quadrics containing the abelian surface and the
scroll $R$. Then we will use a particular version of Lemma
\ref{quadrics} to bound this number.

Let us consider generic hyperplanes $H_1\in \P^{d-1}_1$ and
$H_2\in \P^{d-1}_2$. Let $C_i=\phi_i^*H_i\in |\L_i|$. Note that
$C_1.C_2=2d$. Consider the space $H=\langle H_1,H_2 \rangle$. It
is a $\P^{n-3}$ and it is  invariant by the involution. The base
quadrics of $\P^{n-1}$ intersected with $H$ are base quadrics of
$H$.

Let $F$ be the intersection of $H$ with the scroll $R$. $F$ is
contained in at least $k_b$ independent base quadrics. Because the
lines of $R$ are lines joining the image of the points of
$\bar{A}$ by the maps $\phi_i$, we have:

$$
\bigcup_{P\in C_1\cap C_2} \langle \phi_1(P),\phi_2(P) \rangle\subset F.
$$

The spaces $H_1$ and $H_2$ are generic, so when $d>3$ we can apply
Lemma \ref{tecnico}. We see that $\phi_i(C_1\cap C_2)$ are $2d$
points in general position. The base quadrics of rank $2$ in $H$ are formed by two hyperplanes containing $H_1$ and $H_2$ respectively. It follows that $F$ is not contained in base quadrics of rank $2$.

The space of base quadrics of $H$ has projective dimension
$(d-1)^2-1$ and the subvariety of base quadrics of rank $2$ is of
dimension $2(d-2)$. Thus,

$$
k_b\leq (d-1)^2-1-2(d-2)=d(d-4)+4, \mbox{ when $d>3$.}
$$

When $d=3$ we know that $\phi_i(C_1\cap C_2)$ are points spanning
a line (see Remark \ref{tecnico2}). It follows that $F$ is
contained in at most a finite number of base quadrics of rank $2$.
From this
$$
k_b\leq 2 \mbox{ when $d=3$.}
$$

On the other hand, $k_b=dim(ker(\alpha_b))$. Since
$h^0(\LY)=h^0(\LYY)=d$ and $h^0(\LY\otimes \LYY)=4d$, it follows
that
$$
k_b\geq d^2-4d=d(d-4), \mbox{ when $d>3$.}
$$

Finally, by Proposition \ref{bh} we know that $k_b$ must be
multiple of $d$. Combining the previous bounds with this fact, we
obtain the following possibilities for $k_b$:
$$
\begin{tabular}{|c|c|} \hline
  {Polarization} & {$k_b$} \\ \hline
  {$(1,6)$} & {$0$} \\
  {$(1,8)$} & {$0,4$} \\
  {$(1,10)$} & {$5$} \\
  {$(1,12)$} & {$12$} \\ \hline
\end{tabular}
$$

We can repeat this argument with the three nontrivial torsion
points of $K(\L)$. Let $I_0^+\oplus  I_0^-\oplus I_1^+\oplus
I_1^-$ be the decomposition of the ideal of quadrics containing
$A$ into $\H_{2d}$-modules. Suppose that $d$ is odd. By Proposition \ref{bh} we know:
$$
k_b=dim(I_1^-)+dim(I_0^-)=dim(I_1^-)+dim(I_1^+)=dim(I_1^-)+dim(I_0^+).
$$

Using this fact and the previous bounds for $k$ and $k_b$, one can
compute directly the number of independent quadrics containing the abelian surface $A$ of type $(1,6)$ and $(1,10)$:

\begin{teo}\label{principal1}
Let $A$ be an abelian surface embedded in $\P^{n-1}$ by a very
ample line bundle $\L$ of type $(1,n)$. Then,

\begin{enumerate}

\item If $n=6$, $A$ is not contained in quadrics.

\item If $n=10$, $A$ is contained in $15$ independent quadrics, so
it is projectively normal. \qed

\end{enumerate}
\end{teo}

\subsection{Abelian surfaces of type $(1,8)$ and $(1,12)$.}

Now we will compute the number of independent harmonic quadrics containing the abelian surface $A$. This is equivalent to study the kernel of the map:
$$
\alpha_h:Sym^2(H^0(\LY))\oplus Sym^2(H^0(\LYY))\lrw H^0(\LY^2).
$$
We also consider the maps:
$$
\begin{array}{l}
{\alpha_1:Sym^2(H^0(\LY))\stackrel{\i_1}{\lrw} Sym^2(H^0(\LY))\oplus Sym^2(H^0(\LYY))\stackrel{\alpha_h}{\lrw} H^0(\LY^2)}\\
{\alpha_2:Sym^2(H^0(\LYY))\stackrel{\i_1}{\lrw} Sym^2(H^0(\LY))\oplus Sym^2(H^0(\LYY))\stackrel{\alpha_h}{\lrw} H^0(\LY^2).}\\
\end{array}
$$
We call the images of these maps $Im_1,Im_2$ respectively. 

First we study the polarization of type $(1,8)$. In this case, the abelian
surface $A$ is embedded in $\P^7$, $h^0(\L^2)=32$ and
$h^0(\Te_{P^7}(2))=36$. The possibilities for $k$ are $k=4$ or
$k=8$. The line bundle is projectively normal when $k=4$. The
images of $\bar{A}$ by the maps defined by the line bundles $\LY$
and $\LYY$ are singular octics in $\P^3$ or elliptic quartic
scrolls (see \cite{bilast}). They are not contained in quadrics,
so the maps $\alpha_1$ and $\alpha_2$ are injective. One inmediately sees that $dim(Im_1)=dim(Im_2)=10$ and $h^0(\LY^2)=16$. The image of $\alpha_h$ has at least dimension $10$.

We want to show that $\alpha_h$ is a surjection. Suppose that this is
not true. Then $dim(ker(\alpha_h))>4$. Moreover, since $ker(\alpha_h)$ is an $\H_8$-module, its dimension must be exactly $8$ and the image of $\alpha_h$ is a $12$-dimensional space. This means that the
subspaces $Im_1$ and $Im_2$ intersect in an $8$-dimensional space.

Let us fix a generic divisor $C_2$ in $|\LYY|$. By Lemma
\ref{tecnico}, $\phi_1(C_2)$ is a non-degenerate curve of degree $8$ in $\P^3$ so
it lies at most in a quadric. On the other hand, $\phi_2(C_2)$
lies in a plane, so it is contained in at least $4$ independent
quadrics of rank $2$ in $\P^3$. Thus, there is a $4$-dimensional
subspace $U$ of $Im_2$ corresponding to quadrics containing $\phi_2(C_2)$. If
$dim(Im_1\cap Im_2)=8$, then $dim(U\cap Im_1)\geq 2$. Thus, there
at least two independent quadrics containing $\phi_1(C_2)$, but
this is not possible.

We have checked that the number of harmonic quadrics containing
the abelian surface is exactly $4$. We can repeat this argument
for the three $2$-torsion points in $K(\L)$. Applying Proposition
\ref{bh}, we see that:
$$
4=k_h=dim(I_0^+)+dim(I_1^+)=dim(I_0^+)+dim(I_0^-)=dim(I_0^+)+dim(I_1^-).
$$
It follows that $A$ must be contained in exactly $4$
independent quadrics and we obtain the following result:

\begin{teo}
An abelian surface embedded in $\P^7$ by a very ample line bundle
of type $(1,8)$ is projectively normal. \qed
\end{teo}

Let us study the polarization of type $(1,12)$. Now, the abelian
surface $A$ is embedded in $\P^{11}$, $h^0(\L^2)=48$ and
$h^0(\Te_{P^{11}}(2))=78$. The possibilities for $k$ are $k=30$ or
$k=36$. The line bundle is projectively normal when $k=30$. 

The images of $\bar{A}$ by the maps defined by the line bundles $\LY$
and $\LYY$ are abelian surfaces in $\P^5$ or elliptic normal
scrolls of degree $6$ (see \cite{hula}). 

In the first case, we have seen that they are not contained in quadrics (Theorem \ref{principal1}). Thus, $dim(Im_1)=dim(Im_2)=21$ and $h^0(\LY^2)=24$. It follows that the kernel of $\alpha_h$ has at most dimension $21$. But $ker(\alpha_h)$ is an $\H_{12}$-module, so its dimension is multiple of $6$. We deduce that there are exactly $18$ independent harmonic quadrics containing the abelian surface $A$.

In the second case, it is known that the normal elliptic scrolls in $\P^5$ are contained in $3$ independent quadrics (see \cite{ho1}, \cite{ho2}). From this, $dim(Im_1)=dim(Im_2)=18$ and $h^0(\LY^2)=24$. The kernel of $\alpha_h$ must be of dimension $18$ or $24$. Suppose that it is $24$. Then $Im_1=Im_2$ in $H^0(\LY^2)$. Let us take two generic curves $C,C'\in |\LYY|$. The images $\phi_2(C),\phi_2(C')$ are contained in a degenerate quadric. Because $Im_1=Im_2$, the curves $\phi_1(C),\phi_1(C')$ must be contained in a quadric that does not contain $\phi_1(\bar{A})$. But $\phi_1(\bar{A})$ is an elliptic surface of degree $6$ and by lemma \ref{tecnico}, $\phi_1(C),\phi_1(C')$ are non-degenerate curves of degree $12$, so this is not possible. We conclude that the number of independent harmonic quadrics containing the abelian surface $A$ is again $18$.

Finally, using the number of base quadrics computed in the previous section we see that $A$ is contained exactly in $30$ independent quadrics and we obtain the following result.

\begin{teo}
An abelian surface embedded in $\P^{11}$ by a very ample line bundle
of type $(1,12)$ is projectively normal. \qed
\end{teo}

\end{document}